\documentclass[12pt]{article}
\usepackage{amsmath,amssymb,amsbsy,amsfonts,amsthm,latexsym,
                        amsopn,amstext,amsxtra,euscript,amscd}
\begin{document}
\newtheorem{theorem}{Theorem}
\newtheorem{lemma}[theorem]{Lemma}
\newtheorem{claim}[theorem]{Claim}
\newtheorem{cor}[theorem]{Corollary}
\newtheorem{prop}[theorem]{Proposition}
\newtheorem{definition}{Definition}
\newtheorem{question}[theorem]{Open Question}

%%%%%%%%%%%%%%%%%%%%%%%%%
% Alphabet calligraphic %
%%%%%%%%%%%%%%%%%%%%%%%%%
\def\cA{{\mathcal A}}
\def\cB{{\mathcal B}}
\def\cC{{\mathcal C}}
\def\cD{{\mathcal D}}
\def\cE{{\mathcal E}}
\def\cF{{\mathcal F}}
\def\cG{{\mathcal G}}
\def\cH{{\mathcal H}}
\def\cI{{\mathcal I}}
\def\cJ{{\mathcal J}}
\def\cK{{\mathcal K}}
\def\cL{{\mathcal L}}
\def\cM{{\mathcal M}}
\def\cN{{\mathcal N}}
\def\cO{{\mathcal O}}
\def\cP{{\mathcal P}}
\def\cQ{{\mathcal Q}}
\def\cR{{\mathcal R}}
\def\cS{{\mathcal S}}
\def\cT{{\mathcal T}}
\def\cU{{\mathcal U}}
\def\cV{{\mathcal V}}
\def\cW{{\mathcal W}}
\def\cX{{\mathcal X}}
\def\cY{{\mathcal Y}}
\def\cZ{{\mathcal Z}}

%%%%%%%%%%%%%%%%%%%%%%%
% Alphabet blackboard %
%%%%%%%%%%%%%%%%%%%%%%%
\def\A{{\mathbb A}}
\def\B{{\mathbb B}}
\def\C{{\mathbb C}}
\def\D{{\mathbb D}}
\def\E{{\mathbb E}}
\def\F{{\mathbb F}}
\def\G{{\mathbb G}}
\def\I{{\mathbb I}}
\def\J{{\mathbb J}}
\def\K{{\mathbb K}}
\def\L{{\mathbb L}}
\def\M{{\mathbb M}}
\def\N{{\mathbb N}}
\def\O{{\mathbb O}}
\def\P{{\mathbb P}}
\def\Q{{\mathbb Q}}
\def\R{{\mathbb R}}
\def\S{{\mathbb S}}
\def\T{{\mathbb T}}
\def\U{{\mathbb U}}
\def\V{{\mathbb V}}
\def\W{{\mathbb W}}
\def\X{{\mathbb X}}
\def\Y{{\mathbb Y}}
\def\Z{{\mathbb Z}}

\def\E{{\mathbf E}}
\def\Fp{\F_p}
\def\ep{{\mathbf{e}}_p}
\def\Nm{{\mathrm{Nm}}}
\def\lcm{{\mathrm{lcm}}}

\def\scr{\scriptstyle}
\def\\{\cr}
\def\({\left(}
\def\){\right)}
\def\[{\left[}
\def\]{\right]}
\def\<{\langle}
\def\>{\rangle}
\def\fl#1{\left\lfloor#1\right\rfloor}
\def\rf#1{\left\lceil#1\right\rceil}
\def\le{\leqslant}
\def\ge{\geqslant}
\def\eps{\varepsilon}
\def\mand{\qquad\mbox{and}\qquad}

\newcommand{\comm}[1]{\marginpar{%
\vskip-\baselineskip %raise the marginpar a bit
\raggedright\footnotesize
\itshape\hrule\smallskip#1\par\smallskip\hrule}}

\def\xxx{\vskip5pt\hrule\vskip5pt}

\def\rank{\mathop{rank}}
\def\trace{\mathop{trace}}
\def\ind{\mathrm{ind}}
\def\IM{\mathrm{Im}}

\def\Orb{\mathrm{Orb}(f)}
\def\Orbb{\overline{\mathrm{Orb}}(f)}

%%%%%%%%%%%%%%%%%%
%% PAPER BEGINS %%
%%%%%%%%%%%%%%%%%%

\title{{\bf  On the Length of Critical Orbits of Stable 
Quadratic Polynomials}}

\author{ 
{\sc Alina~Ostafe}\\
{Institut f\"ur Mathematik, Universit\"at Z\"urich}\\
{Winterthurerstrasse 190 CH-8057, Z\"urich, Switzerland}\\
{\tt alina.ostafe@math.uzh.ch}
\and
{\sc Igor E.~Shparlinski} \\
{Department of Computing, Macquarie University} \\
{Sydney, NSW 2109, Australia} \\
{\tt igor@ics.mq.edu.au}}

\date{\today}
\pagenumbering{arabic}

\maketitle

\begin{abstract} We use the Weil bound of multiplicative character sums
together with some recent results of N.~Boston and  R.~Jones, 
to show that the   {\it critical orbit\/}
of quadratic polynomials over a  finite field of $q$ 
elements 
is of length    $O\(q^{3/4}\)$, 
improving upon the trivial bound $q$. \end{abstract}

\section{Introduction}  

Let $\F_q$ be a finite field of $q$ elements. For a polynomial 
$f\in \F_q[X]$ we define the sequence of iterations:
$$
f^{(0)}(X)  = X, \qquad f^{(n)}(X)  = f\(f^{(n-1)}(X)\), \quad 
n =1, 2, \ldots\,.
$$
Following~\cite{Ali,AyMcQ,Jon2,JB}, we say that $f$ is {\it stable\/} if all
polynomials $f^{(n)}$ are irreducible over $\F_q$. 

We now assume that $q$ is odd.

As in~\cite{JB}, for a quadratic 
polynomial $f(X) = aX^2 + bX + c \in \F_q[X]$, $a \ne 0$, 
we define $\gamma= -b/2a$ as the unique critical 
point of $f$ (that is, the zero of the derivative $f'$) 
and consider the set
$$
\Orb = \{f^{(n)}(\gamma)\ : \ n =2, 3,   \ldots\}
$$
which is called the {\it critical orbit\/} of $f$. 
Clearly there is some $t$ such that $f^{(t)}(\gamma) = f^{(s)}(\gamma)$
for some positive integer $s < t$. 
Then   $f^{(n+t)}(\gamma) = f^{(n+s)}(\gamma)$ for any $n \ge 0$.
Accordingly, for the smallest value of $t_f$ with the above
condition,  we have 
$$
\Orb = \{f^{(n)}(\gamma)\ : \ n = 2,   \ldots, t_f\}
$$
and $\# \Orb = t_f-1$ or $\# \Orb = t_f-2$ (depending whether $s =1$ or 
$s \ge 2$ in the above).
It is shown in~\cite{Jon1,Jon2,JB} that critical orbits
play a very important role in the dynamics
of polynomial iterations. 

Trivially we have $t_f \le q+1$. In fact, by the
{\it Birthday Paradox\/} one expects that $t_f$ is
of order $q^{1/2}$ and there are examples of 
polynomials which have orbits of about this length. 

Here we obtain a nontrivial upper bound on the orbit length of 
stable quadratic polynomials:

\begin{theorem}
\label{thm:UB} For any odd $q$ and  any stable quadratic polynomial 
$f \in\F_q[X]$ we have
$$
t_f = O\(q^{3/4}\).
$$
\end{theorem}

By~\cite[Proposition~3]{JB}, a quadratic polynomial $f\in\F_q[X]$  
is stable if the {\it adjusted orbit}
$$
\Orbb = \{-f(\gamma)\} \bigcup \Orb
$$
contains no squares. 
We also recall that 
$\alpha \in \F_q$ is a square if either $\alpha = 0$ or 
$\alpha^{(q-1)/2}=1$ that can be tested (via
repeated squaring) in $O(\log q)$ field operations. 
Combining these with the bound of Theorem~\ref{thm:UB}, 
 we immediately obtain:

\begin{cor}
\label{cor:Test} For any odd $q$, a quadratic polynomial 
$f\in \F_q[X]$ can be tested for stability in time 
$q^{3/4+o(1)}$. 
\end{cor}

Our proof is based on  the Weil bound for character sums  
with polynomials, see~\cite[Theorem~11.23]{IwKow}.

Finally, we  remark that estimating the size 
of the set of stable quadratic polynomials 
$aX^2 + bX + c \in \F_q[X]$ is a very interesting question
to which we hope our technique can apply as well.

\section{Proof of Theorem~\ref{thm:UB}} 

Let $\chi$ be the quadratic character of $\F_q$,

By~\cite[Proposition~3]{JB}, if a  quadratic polynomial
$f \in\F_q[X]$ is   stable then $\Orb$  contains no squares, 
that is, $\chi\( f^{(n)}(\gamma)\) = -1$, $n=2,3,\ldots$.

We now fix an integer parameter $K$ and note that for any $n\ge 1$,
we have simultaneously
$$
\chi\( f^{(k+n)}(\gamma)\) = -1, \qquad k =1,\ldots, K,
$$
which we rewrite as 
\begin{equation}
\label{eq:K cond}
\chi\( f^{(k)}\( f^{(n)}(\gamma)\)\) = -1, \qquad k =1,\ldots, K.
\end{equation}
Since by the definition of $t_f$, the values $f^{(n)}(\gamma)$, $n=1, \ldots, t_f-1$,
are pairwise distinct elements of $\F_q$ we derive from~\eqref{eq:K cond} that
\begin{equation}
\label{eq:t and T}
t_f-1 \le \# \cT_q(K)
\end{equation}
where 
$$
\cT_q(K) =  \left\{x \in   \F_q~:~ \chi\(f^{(k)}(x)\) = -1, 
\ k=1, \ldots, K\right\}.
$$
We have
\begin{equation}
\label{eq:TqK}
\# \cT_q(K) = \frac{1}{2^K} \sum_{x \in   \F_q} \prod_{k=1}^K \(1-\chi\(f^{(k)}(x)\)\)
\end{equation}
since for every $x \in  \cT_q(K)$ the product 
on the right hand side of~\eqref{eq:TqK}  is $2^K$; otherwise
it is  $0$ when $\chi(f^{(k)}(x)) =1$ for at least 
one $k =1, \ldots, K$ (note that since by our assumption $f^{(k)}(X)$ is 
irreducible over $\F_q$ we have $f^{(k)}(x) \ne 0$ for $x \in   \F_q$).

Just expanding the product in~\eqref{eq:TqK} and 
changing the order of summation, we obtain $2^k-1$  character sums 
 of the shape
\begin{equation}
\label{eq:Sums}
(-1)^\nu \sum_{x \in   \F_q}  \chi\(\prod_{j=1}^\nu f^{(k_\nu)}(x)\),
\qquad 1 \le k_1 < \ldots < k_\nu\le K,
\end{equation}
with $\nu \ge 1$ and one trivial sum which equal to $q$ (corresponding to
the terms $1$ in the product in~\eqref{eq:TqK}).

Clearly $f^{(k)}(X)$ is a polynomial of  degree
$2^k$. By our assumption they are irreducible, therefore none of the polynomials
$$
\prod_{j=1}^\nu f^{(k_\nu)}(X) \in \F_q[X],
\qquad 1 \le k_1 < \ldots < k_\nu\le K,
$$ 
is a perfect square. Therefore the {\it Weil bound\/} 
see~\cite[Theorem~11.23]{IwKow}, applies to every 
sum~\eqref{eq:Sums} and implies that each of them is $O(2^K q^{1/2})$.
Therefore
\begin{equation}
\label{eq:T bound}
\# \cT_q(K) = \frac{1}{2^K} q + O(2^K q^{1/2}).
\end{equation}
Choosing $K$ to satisfy
$$
2^K \le q^{1/4} < 2^{K+1}
$$
and combining~\eqref{eq:t and T} and~\eqref{eq:T bound}
we conclude the proof.

\section{Comments}

It is certainly interesting to obtain nontrivial estimates 
on the size $S_q$ of the set  of the triples 
$(a,b,c) \in \F_q^* \times \F_q   \times \F_q$ which 
correspond to stable quadratic polynomials $f(X) = aX^2 + bX + c$.
Denoting by $F_k(a,b,c)$ the $k$th element of the critical 
orbit of $f$, we see that for any integer parameter $K$ we have
\begin{equation}
\label{eq:S and W}
S_q\le \# \cW_q(K),
\end{equation}
where 
$$\cW_q(K) = \left\{(a,b,c) \in \F_q^* \times \F_q  \times \F_q~:~ \chi\(F_k(a,b,c)\) = -1, 
\ k=1, \ldots, K\right\}, 
$$
and  as before $\chi$ denotes the quadratic character of $\F_q$.
As in the proof of Theorem~\ref{thm:UB}, we have
\begin{equation}
\label{eq:WqK}
\# \cW_q(K) \le \frac{1}{2^K} \sum_{(a,b,c) \in \F_q^* \times \F_q  \times \F_q} \prod_{k=1}^K \(1-\chi\(F_k(a,b,c)\)\)
\end{equation}
since for every triple $(a,b,c) \in  \cW_q(K)$ the product 
on the right hand side of~\eqref{eq:WqK}  is $2^K$; otherwise
it is either $0$  (when $\chi(F_k(a,b,c)) =1$ for at least 
one $k =1, \ldots, K$) or $1$ (when $F_1(a,b,c) =\ldots = F_K(a,b,c) = 0$).

Clearly $F_k(a,b,c)$ are rational functions in $a,b,c$ of degree at most 
$O(2^k)$. Just expanding the product in~\eqref{eq:WqK} and 
changing the order of summation, we obtain $2^k-1$  character sums 
 of the shape
\begin{equation}
\label{eq:Sums2}
(-1)^\nu \sum_{(a,b,c) \in \F_q^* \times \F_q  \times \F_q}  \chi\(\prod_{j=1}^\nu F_{k_\nu}(a,b,c)\),
\quad 1 \le k_1 < \ldots < k_\nu\le K,
\end{equation}
with $\nu \ge 1$ and one trivial sum corresponding to  $1$
in~\eqref{eq:WqK}. 
Assuming that one can prove that the Weil-type bound 
$O(2^K q^{5/2})$ applies to 
all of them, we obtain from~\eqref{eq:S and W} that 
$S_q = O(q^3/2^K + 2^K q^{5/2})$ and optimising the
choice of $K$ we derive $S_q =O(q^{11/4})$. In fact, for a nontrivial
estimate of $S_q$ it is enough to show that almost all sums~\eqref{eq:Sums2}
admit a  nontrivial estimate, which we pose as an open question. 

\section*{Acknowledgement}

The authors are grateful to Rafe Jones and 
Arne Winterhof for careful reading of the
preliminary version of the 
manuscript and many useful comments.

During the preparation of this paper,  
A.~O. was supported in part by 
the Swiss National Science Foundation   Grant~121874  
and I.~S. by
the  Australian Research Council 
Grant~DP0556431.

\end{document}